\documentclass[12pt]{article}
\usepackage{graphicx,paralist,amsmath,amssymb,amsthm}
\newtheorem{thm}{Theorem}[section]
\newtheorem{cor}[thm]{Corollary}
\newtheorem{lem}[thm]{Lemma}

\theoremstyle{definition}
\newtheorem{defn}[thm]{Definition}
\newtheorem{rem}[thm]{Remark}
\numberwithin{equation}{section}
\newcommand{\abs}[1]{\left\vert#1\right\vert}
\newcommand{\A}{\mathcal{A}}
\newcommand{\F}{\mathcal{F}}
\newcommand{\K}{\mathcal{K}}
\newcommand{\N}{\mathcal{N}}
\newcommand{\M}{\mathcal{M}}
\newcommand{\LL}{\mathcal{L}}
\newcommand{\PP}{\mathcal{P}}
\newcommand{\HH}{\mathcal{H}}
\newcommand{\RR}{\mathcal{R}}

\begin{document}

\title{Maximum Distance Separable Codes\\ and Arcs in Projective Spaces}%
\author{T. L. Alderson \thanks{The author acknowledges support from the N.S.E.R.C. of Canada} \\
Mathematical Sciences\\
University of New Brunswick\\
Saint John, NB.\\
E2L 4L5\\
Canada
 \and
A. A. Bruen\thanks{The author acknowledges support from the N.S.E.R.C. of Canada}\\
Mathematics and Statistics\\
University of Calgary\\
Calgary, AB\\
T2N 1N4\\
Canada
 \and
R. Silverman\\
Mathematics\\
Wright State University\\
Dayton, OH.\\
45431\\
U.S.A.
}%

\date{}

\begin{abstract}

\noindent Given any linear code $C$ over a finite field $GF(q)$ we
show how $C$ can be described in a  transparent and geometrical
way by using the associated Bruen-Silverman code.

Then, specializing to the case of MDS codes we use our new approach
to offer improvements to the main results currently available
concerning MDS extensions of linear MDS codes. We also sharply limit
the possibilities for constructing long non-linear MDS codes. Our
proofs make use of the connection between the work of R\'{e}dei
\cite{Rdi} and the R\'{e}dei blocking sets that was first pointed
out over thirty years ago in  \cite{BrL}. The main results of this
paper significantly strengthen those in \cite{BBT},\cite{BTB}.
\end{abstract}
 \maketitle

\section{Introduction}

Maximum distance separable codes (MDS codes) are at the heart of
combinatorics and finite geometries.  In their book \cite{M-S}
MacWilliams and Sloane describe MDS codes as ``one of the most
fascinating chapters in all of coding theory".  These codes can be
linear or non-linear.  We define them as follows.

\begin{defn}
An $(n,k,q)$-MDS code $C$ is a collection of $q^k$  $n$-tuples (or
code words) over an alphabet $\A$ of size $q$ satisfying the
following condition:  No two words of $C$ agree in as many as $k$
coordinate positions.
\end{defn}

In the special case that $\A=GF(q)$ the finite field of order $q$
and $C$ is a linear vector space of dimension $k$ and length $n$, we
say that $C$ is a \textit{linear} $(n,k,q)$-MDS code.  These linear
MDS codes are fundamental in error correction.  In part this is due
to the fact that they are precisely the linear codes meeting the
Singleton bound (see \cite{BrF}) which states that the minimum
Hamming distance $d$ of a linear $(n,k,q)$-MDS code satisfies
$d=n-k+1$. In particular, a class of these MDS codes, the so-called
Reed-Solomon codes, are a mainstay of modern industrial
applications. Geometrically, these Reed-Solomon codes are precisely
the normal rational curves. Fundamental new algorithmic results on
decoding Reed-Solomon codes are described in the work of Madhu
Sudan, the winner of the 2002 Nevanlinna Prize (see \cite{Wig}).
Concerning the possible applications of error correction we should
also mention the important recent results of Calderbank and Shor
\cite{Cal} relating to the emerging area of quantum error correcting
codes.

The following combinatorial result is shown in \cite{Sil}.

\begin{lem}\label{Lsil}
Let $C$ be an $(n,k,q)$-MDS code over the alphabet $\A$.
\begin{enumerate}[(a)]

\item \label{pt1} Fix any $t$ coordinate positions
$a_1<a_2<\cdots<a_t$, $t\leq k$, and choose $\alpha_1,
\alpha_2,\ldots,\alpha_t\in \A$ (not necessarily distinct). Then
there are exactly $q^{k-t}$ code words in $C$ such that the entry
in position $a_i$ equals $\alpha_i$, $1\leq i \leq t$.

\item \label{pt2}Fix $\alpha \in \A$ and fix some coordinate
position $j$. Let $C_1$ denote the set of all code words in $C$
having $\alpha$ in position $j$. Then by deleting the $j$'th
coordinate position from $C_1$ we obtain an $(n-1,k-1,q)$-MDS
code.

\item \label{pt3}If $q$ is even then $n\leq q+k-1$, if $q$ is odd
and $k>2$ then $n\leq q+k-2$.
\end{enumerate}
\end{lem}

The existence question for codes meeting the combinatorial bound in
Lemma \ref{Lsil} (\ref{pt3}) is largely an open one. For $k=2$ we
have $n\leq q+1$, with $n=q+1$ if and only if there exists an affine
plane of order $q$.  Such planes exist if $q$ is a prime power.
Whether or not this condition is necessary (i.e. `` the prime power
conjecture") has been one of the most important and outstanding open
questions in finite geometries for over 50 years, since the
publication of the Bruck-Ryser theorem.  For $k=3$, equality is also
possible.

For $k=4$ the only known result is that in the case of equality,
36 divides $q$ \cite{B&S3}, (so that $q$ could not be a prime
power in this case).  From the results of \cite{Ald1} it follows
that for $k\geq 4$, if $q$ is even and 36 does not divide $q$ then
$n \leq q+k-3$.

In searching for ``long" MDS codes a natural approach is to begin
with a fixed code $C$ and attempt to ``lengthen" the code, while
preserving the MDS property.

\begin{defn} Let $C$ be an $(n,k,q)$-MDS code.
A code  $\,C'$ is said to be an \textit{extension} of $\,C$ if
$C'$ is an $(m,k,q)$-MDS code where $m>n$ and upon deleting some
fixed set of $m-n$ coordinate positions from each word of $\,C'$
the code $\,C$ is obtained. Equivalently, $C$ is said to be
\textit{extendable} (to the code $C'$).  An MDS code is said to be
\textit{maximal} if it admits no extensions.
\end{defn}

The following is an immediate consequence of Lemma \ref{Lsil}
\textit{(b)}.

\begin{lem}\label{Lptx}
An $(n,k,q)$-MDS code $C$ is extendable to an $(n+1,k,q)$-MDS code
if and only if there exists a  partition
$\PP=\{C_1,C_2,\ldots,C_q\}$ of $C$ such that each $C_i$ is an
$(n,k-1,q)$-MDS code.
\end{lem}

\begin{defn}
An \textit{$n$-arc} $\K$ in $PG(k,q)$ is a collection of $n \geq
k+1$ points no $k+1$ of which are incident with a common hyperplane.
A \textit{dual $n$-arc} in $PG(k,q)$ is a collection of $n \geq k+1$
hyperplanes  no $k+1$ of which are incident with a common point.
\end{defn}

Suppose $C$ is a linear $(n,k,q)$-MDS code (so $\A=GF(q)$) and
choose a generator matrix $G$ for $C$ of size $k\times n$.  The MDS
property of $C$ is equivalent to the condition that every collection
of $k$ columns of $G$ is linearly independent (see \cite{BrF}).  As
observed by B. Segre one can multiply the columns of $G$ by nonzero
scalars and still preserve the MDS property. The columns of $G$ can
therefore be regarded as points (or, by duality, as hyperplanes)
belonging to an $n$-arc in $PG(k-1,q)$. \textit{Hence any results on
linear MDS codes can be translated to an equivalent theorem on
arcs.}  A \textit{normal rational curve} in $PG(k,q)$, $2\leq k \leq
q-2$ is a $(q+1)$-arc projectively equivalent to the $(q+1)$-arc
$\{(1,t,\ldots,t^k)\,|\, t\in GF(q)\,\} \cup\{(0,\ldots,0,1)\}$. The
$n$-arcs which are subsets of normal rational curves correspond to
\textit{generalized Reed-Solomon} (GRS) codes. Linear
$(q+1,k,q)$-MDS codes are therefore easily constructed.

Denote by $m(k,q)$ the size of the largest (dual) arc in
$PG(k,q)$. Finding the value of $m(k,q)$ has been the focus of
much research (see \cite{BBT,BTB,Cas,JAT,Segr}).  The Main
Conjecture for linear MDS codes, always taking $q>k$, is the
following:
\[ m(k-1,q)=\left\{\begin{array}{ll}
                    q+2 & \text{ if  $k=3$ and $k=q-1$ both with $q$ even} \\
                    q+1 & \text{ in all other cases}
                  \end{array}\right.
                \]

The Main Conjecture has not been proved in general.  It has been
verified in many cases. In their paper \cite{BTB} Bruen, Thas, and
Blokhuis show it to hold at least asymptotically. The existence and
possible structure of long MDS codes was a central theme in the
address of J.A. Thas to the International Congress of Mathematicians
in 1998 \cite{JAT}.  The question of the existence of MDS codes
meeting the combinatorial bound remains largely an open one.  In his
engineering textbook, (\cite{Wick}, page 193), Wicker describes the
analogous question for arcs as ``one of the more interesting
problems in projective geometry over Galois fields".

Consider a linear $(n,k,q)$-MDS code $C$ over $\F=GF(q)$ with
associated generator matrix $G$.  A linear extension of $C$ arises
by augmenting $G$ with an appropriate column vector. Over $\F$
there are in total $q^k$ column vectors to check using (perhaps
naively) an  exhaustive search.  In order to consider general (not
necessarily linear) extensions of $C$ we let $M$ be a $q^k \times
n$ array whose rows are precisely the words of $C$. This is called
the \textit{matrix form} of $C$.  A general extension of $C$
arises by augmenting $M$ with an appropriate column vector. Over
$\F$ there are a total of $q^{q^k}$ possible column vectors.
Hence, the search for an extension of $C$ grows exponentially when
one considers general as well as linear extensions.

Our main result (Theorem \ref{C8b}) shows that linear MDS codes of
reasonable length do not admit non-linear extensions.  In
particular this gives a strengthening of results in the two papers
mentioned above.  Consequently, in the search for MDS codes close
to the combinatorial bound, the approach of extending long linear
codes in a non-linear fashion is  fruitless. Moreover, by results
such as those of Bruen, Thas and Blokhuis \cite{BTB}, there are
many linear MDS codes of reasonable length that admit no linear
extensions.  Now we can show that these codes admit no extensions
whatsoever.

\section{2-Dimensional MDS Codes, Bruck Nets}

A \textit{Bruck net} (see \cite{BrF}) $\N$ is a finite incidence
structure of points and certain subsets of points called lines
satisfying the following axioms.
\begin{enumerate}[(a)]
\item Any two points are contained in (lie on) at most one line.
 \item Given a line $\ell$ and point $P$ off $\ell$, there is a
 unique line through $P$ failing to meet $\ell$.
 \item There exists a triangle in $\N$ i.e. a set of three
 non-collinear points $A,B,C$ such that the point pairs $AB$, $AC$,
 and $BC$ are joined in $\N$.
\end{enumerate}

Let us suppose that some line of $\N$ contains exactly $q$ points.
Then all lines of $\N$ contain exactly $q$ points and $q$ is called
the \textit{order} of $\N$. From axiom (b) it follows that the lines
of $\N$ fall into $n$ \textit{parallel classes}.  Each parallel
class has exactly $q$ lines no two of which meet.  Two lines from
distinct parallel classes meet in a unique point. The total number
of points in $\N$ is $q^2$ and the total number of lines is $nq$.
The parameter $n$ is called the \textit{degree} of $\N$. To extend a
net $\N$ is to append one or more parallel classes of lines (thereby
increasing the degree). It can be shown that $n\leq q+1$ with
equality if and only if $\N$ consists of the points and lines of an
affine plane of order $q$. In particular, a net of order $q$ and
degree $q+1$ can not be extended.

\begin{lem}
A Bruck net of order $q$ and degree $n$ is equivalent to a
$(n,2,q)$-MDS code.
\end{lem}

This is shown in \cite{BrF}.  Briefly, each of the $q^2$ points
$P$ of $\N$ give a code word as follows.  Label the $n$ parallel
classes of $\N$ as $\{1,2,\ldots,n\}$ where $n\geq3$. Within each
parallel class, label the lines from $\{1,2,\ldots,q\}$. If the
point $P$ lies on the line $\alpha_i$ from parallel class number
$i$, $1\leq i \leq n$, we associate with $P$ the code word
$(\alpha_1,\alpha_2,\ldots,\alpha_n)$.  Two points are joined in
the net if and only if the corresponding code words share a common
entry.

Each coordinate of an $(n,2,q)$-MDS code corresponds to a parallel
class of lines in the related net and so we have a natural
correspondence between extending the code and extending the net.

\section{Linear Codes, BRS Codes}\label{LinBRS}

In the sequel we need to discuss equivalence of codes. Let $C_1$ and
$C_2$ be codes of length $n$ over an alphabet $\A$.  Identify each
code with a matrix, the rows of each matrix being the code words.
Then $C_2$ is said to be \textit{equivalent} to $C_1$ if $C_2$ can
be obtained from $C_1$ by a sequence of operations of the following
three types:
\begin{enumerate}
\item A permutation of the rows of $C_1$\,;
 \item A permutation of the columns of $C_1$\,;
 \item A permutation of the alphabet $\A$ is applied (entry-wise) to a column of $C$.
\end{enumerate}

If two codes are equivalent then corresponding parameters such as
minimum distance are equal, so the codes are essentially
identical. A code that is equivalent to a linear code is referred
to as a code of \textit{type LE}, or an \textit{LE code}. An LE
code need not be linear. For example, if we suitably permute the
symbols in a given column of a linear code, the resulting code
will not contain the zero vector and will therefore not be linear.

Let $C$ be any linear code of length $n$ and dimension $k$ over
the finite field $\F=GF(q)$.  Associated with $C$ is a $k\times n$
generator matrix $G$ of rank $k$ over $\F$.  Each vector of $C$ is
a linear combination of the rows of $G$.  Denote the entries of
$G$ as follows.

\[  G=\left[
\begin{array}{cccc}
  a_{11} & a_{12} & \cdots & a_{1n} \\
  a_{21} & a_{22} & \cdots & a_{2n} \\
  \vdots & \vdots & \ddots & \vdots \\
  a_{k1} & a_{k2} & \cdots & a_{kn} \\
\end{array}%
\right]
\]
Then a code word $w$ of $C$ can be written as
\begin{equation}\label{eq1}
w=\sum_{i=1}^k \alpha_iR_i
\end{equation}
where $R_i$ denotes the $i^{th}$ row of $G$.\\

We want to get a better geometrical picture of $C$.  This can be
done as follows.\\

Associate with $C$ a projective space $\Sigma=PG(k,q)$ having
homogeneous coordinates $(x_1,x_2,\ldots,x_{k+1})$.  Assume the
hyperplane at infinity $\Sigma_{\infty}$ has equation $x_{k+1}=0$.
So $\Sigma_{\infty}$ has projective dimension $k-1$.  Each column in
$G$, say the $i^{th}$ column, gives rise to a hyperplane $\Pi_i$ in
$\Sigma_{\infty}$.  The subspace $\Pi_i$ of projective dimension
$k-2$ is defined to be the solution set of the following system of
equations:
\[ \left\{ \begin{array}{l}
             x_{k+1}=0, \\
             a_{1i}x_1+a_{2i}x_2+\cdots a_{ki}x_k=0. \\
           \end{array} \right.\]

Let $E=\Sigma\setminus\Sigma_{\infty}$ denote the associated
affine (or vector) space of dimension $k$.  Thus $E$ has $q^k$
points or vectors.  Each point $P$ in $E$ has homogeneous
coordinates $(\alpha_1,\alpha_2,\ldots,\alpha_{k}, 1)$.  We wish
to associate with $P$ a code word
$(\lambda_1,\lambda_2,\ldots,\lambda_n)$.  We have that $P$ lies
on a certain hyperplane labelled $H_i(P)$ containing the subspace
$\Pi_i$ for each $i$, $1\leq i \leq n$. If we label the $q$
hyperplanes of $\Sigma $ other than $\Sigma_\infty $ containing
$\Pi_i$, then $P$ will lie on say the hyperplane labelled
$\lambda_i$.  In this way we will end up with a code $C_1$
consisting of $q^k$ code words
$(\lambda_1,\lambda_2,\ldots,\lambda_n)$ of length $n$ over $\F$.
$C_1$ will be a \textit{Bruen-Silverman code} (BRS code)
associated with $C$.

Note that the construction of $C_1$ is analogous to the
construction of the MDS code associated with a Bruck net in the
previous section.  The idea of  BRS type codes was first
introduced in \cite{B&S3}.

The code $C_1$ will depend on the labelling of $H_i(P)$.  To fix
coordinates, let us proceed in the following way.  Assume that the
point $U=(0,1,0,\ldots,0)$ is not contained in any of the
subspaces $\Pi_i$, $1\leq i \leq n$.  It follows that in $G$,
$a_{2i}$ is nonzero $1\leq i \leq n$. Multiplying any column of
$G$ by a nonzero constant yields a generator matrix for a code
equivalent to $C$. Hence, we may assume $a_{2i}=1$, $1\leq i \leq
n$. Let $V$ denote the point $(0,0,\ldots,0,1)$ so that $V$ is a
point of $E$. Then any point of $E$ on the line $\ell=UV$ has
coordinates $(0,\gamma,0,\ldots,0,1)$.  Let
$P=(\alpha_1,\alpha_2,\ldots,\alpha_k,1)$ be an affine point.  The
hyperplane $H_i(P)$ containing $P$ and $\Pi_i$ meets $\ell$ in a
point $(0,\gamma_i,0,\ldots,0,1)$ and $\gamma_i$ gives the
$i^{th}$ coordinate entry of the code word
$(\gamma_1,\gamma_2,\ldots,\gamma_n)$ in $C$ associated with the
point $P$.

Let us calculate $\gamma_i$.  Any hyperplane other than
$\Sigma_{\infty}$ containing $\Pi_i$ has an equation of the form

\begin{equation}\label{eq2}
a_{1i}x_1+ a_{2i}x_2+\cdots + a_{ki}x_k + b x_{k+1}=0
\end{equation}

If the hyperplane contains
$P=(\alpha_1,\alpha_2,\ldots,\alpha_k,1)$ then we have
\[a_{1i}\alpha_1+ a_{2i}\alpha_2+\cdots + a_{ki}\alpha_k +
b=0\] which gives
\begin{equation}\label{eq3}
b= -\left( a_{1i}\alpha_1+ a_{2i}\alpha_2+\cdots +
a_{ki}\alpha_k\right)
\end{equation}

The hyperplane given by (\ref{eq2}) meets the line $\ell$ in the
point $(0,\gamma_i,0,\ldots,0,1)$ so that $\gamma_ia_{2i}+b=0$.
Since we have $a_{2i}=1$ we get
\[\gamma_i=-b=a_{1i}\alpha_1+ a_{2i}\alpha_2+\cdots +
a_{ki}\alpha_k\]  But then the code word of $C_1$ associated with
$P$ is precisely the code word of $C$ in (\ref{eq1}) associated
with the coefficients $\alpha_1,\alpha_2,\ldots,\alpha_k$.  We
have shown the following.

\begin{thm}
The code $C_1$ is identical to the original code $C$.  In particular
$C_1$ is linear.
\end{thm}

To summarize, we now have a completely new way of looking at the
original code $C$, as follows.  We can identify a code word $w$ in
$C$ with the set of coefficients
$\alpha_1,\alpha_2,\ldots,\alpha_k$ as in formula \ref{eq1}.
Alternatively, the code word can be thought of as a point
$P=(\alpha_1,\alpha_2,\ldots,\alpha_k,1)$ in an affine space.  To
find the $i^{th}$ coordinate of $w$, given $P$, we calculate the
label of the unique hyperplane $H_i$ of $\Sigma$ (the underlying
projective space) containing $\Pi_i$ and $P$ (using $\ell$) as
above. Here $\Pi_i$ is a subspace of $\Sigma$ of codimension 2
corresponding to the $i^{th}$ column of $G$, the generator matrix
of $C$.\\

The power of this new approach will be demonstrated in section
\ref{l3d}.  From this picture it is clear that the set of code
words with a given symbol in the $i^{th}$ coordinate position
corresponds to the points of $E$ contained in a certain hyperplane
$H_i$ of $\Sigma$. The code words with given symbols in two fixed
positions $i$ and $j$ correspond to the intersection $H_i\cap H_j$
of two hyperplanes, and so on.  Hence, two code words $w_1$ and
$w_2$ corresponding to the affine points $P$ and $Q$ will have $t$
common entries if and only if the line $PQ$ intersects
$\Sigma_{\infty}$ in a point belonging to $t$ of the $\Pi_i$'s.\\

In the sequel we will need the following concept.

\begin{defn} \label{dtrs} Let $\K$ be a dual arc
in $\Pi=PG(k,q)$, $k\geq 2$ and let $ \Sigma=PG(k+1,q)$.  A point
set $S$ in $\Sigma-\Pi$ is called a \textit{transversal set} of
$\K$ if every line of $\Sigma$ on a $k$-fold point of $\K$
intersects $S$ in at most one point.  Here, a $k$-fold point of
$\K$ is a point incident with precisely $k$ members of $\K$.
\end{defn}

In the special case that the columns of $G$ correspond to a dual
arc we see that a transversal set corresponds to a collection of
code words, no two of which agree in as many as $k$ coordinates.

\section{Primitive Sets, Slope Sets, Directions.}

Let $\pi$ be any projective plane of order $q$ say and let $\ell$ be
a line of $\pi$.  Let $S$ be a set of points of $\pi$ not on $\ell$,
so $S\subset \pi\setminus\{\ell\}$.  The \textit{R\'{e}dei set} of
$S$ with respect to $\ell$, denoted by $R_{\ell}(S)$, is defined to
be the set of all points of the form $PQ\cap\ell$, where $P $ and
$Q$ are distinct points of $S$ and $PQ$ denotes the line joining
them.

\begin{defn}\label{Dprim} Let $\pi$, $\ell$ and $S$ be as above with
$|S|=q$. Let $A$ be a subset of the points of $\ell$. Then $A$ is
said to be \textit{primitive} if whenever $R_{\ell}(S)\subset A$,
the set $S$ must be contained in a line of $\pi$.
\end{defn}

\begin{thm}\label{T2}
Let $\pi$ be a projective plane of order $q$ (not necessarily a
prime power) with a distinguished line $l_{\infty}$. Let $A$ be a
subset of $l_{\infty}$ and $S$ a set of $q$ points of
$\pi\setminus \{l_{\infty}\}$ with the following property. The
line joining any two points of $S$ intersects $l_{\infty}$ within
$A$ i.e. $R_{\ell}(S)\subset A$ . If $|A| <\sqrt{q}+1$, then $S$
is a subset of a line of $\pi$.
\end{thm}
\begin{proof}
See: \cite{Br1}, \cite{Br2} and \cite{BrL}. The result is also
implicit in earlier work of R.H. Bruck and T.G. Ostrom.
\end{proof}

\begin{cor}
Let $\pi$ be a projective plane of order $q$ with a distinguished
line $l_{\infty}$. Let $A$ be a subset of $l_{\infty}$. If
$|A|<\sqrt{q}+1$, then $A$ is primitive.
\end{cor}

\begin{thm}\label{T3}
Let $\pi=PG(2,p)$, $p$ a prime with a distinguished line $\ell$ at
infinity. Let $S$ be a set of $p$ affine points and let $A\subset
\ell$ with $R_{\ell}(S)\subset A$. Then if $|A|<(p+3)/2$, S is a
subset of a line of $\pi$.
\end{thm}
\begin{proof}
See: \cite{LOV}.
\end{proof}

\begin{defn} For $q$ a prime power we denote by $\PP(q)$ the size of
the smallest non-primitive set of (collinear) points in $PG(2,q)$.
\end{defn}

\begin{thm}\label{T3b}
Let $\pi=PG(2,q)$ with a distinguished line $l_{\infty}$. Let
$q=p^h$, $p$ a prime and let $t<h$ be maximal such that $t$ divides
$h$. Let $S$ be a set of $q$ affine points and let $A\subset \ell$
with $R_{\ell}(S)\subset A$. Then if $|A|<p^{t-h}+1$, S is a subset
of a line of $\pi$.
\end{thm}
\begin{proof}
See: \cite{BBBSS},\cite{Ball}.
\end{proof}

\noindent From Theorems \ref{T3} and \ref{T3b} we get the following.

\begin{cor} \label{smallest non-primitive set} Let $\pi=PG(2,q)$ where
$q=p^h$, $p$ prime, and let $t<h$ be maximal such that $t$ divides
$h$. Let $A$ be a set of points on a line $\ell$. If $|A|<\epsilon$
where
\begin{equation} \epsilon=\left\{\begin{array}{ll}
 \frac{q+3}{2} &  h=1\\
p^{h-t}+1 & \text{ otherwise.}
\end{array}\right.\end{equation}
then $A$ is primitive.  In particular we have
\begin{equation*}
\PP(q)\geq \left\{\begin{array}{ll}
 \frac{q+3}{2} &  h=1\\
p^{h-t}+1 & \text{ otherwise.}
\end{array}\right.\end{equation*}
\end{cor}

Let $\pi=PG(2,q)$ with a distinguished line $\ell$ and let $S$ be a
set of $q$ points in $\pi\setminus\{\ell\}$.  If the points of $S$
do not determine all directions (i.e. if $R_{\ell}(S)\subsetneq
\ell$) then $S$ may be regarded as a function on $GF(q)$. (Briefly,
let $\LL$ be the set of all lines incident with at least two points
of $S$. By assumption there exists a point $Q\in \ell$ not incident
with any line of $\LL$. Assign homogeneous coordinates
$(x_1,x_2,x_3)$ in such a way that $\ell$ is given by $x_3=0$ and
$Q=(0,1,0)$. Then $S=\{(a_i,b_i,1)\,|\, i=1\ldots q\}$ yields a
function $f(x)$ where $f(a_i)=b_i$.) Conversely, any function can be
regarded as such a set $S$. Associated with a function $f$ is the
corresponding R\'{e}dei set $\RR$ where
\begin{equation}
\RR=\left.\left\{\frac{f(a_i)-f(a_j)}{a_i-a_j}\,\right|\, i \neq
j\right\}
\end{equation}

\noindent Thus, according to the Definition \ref{Dprim} we have
the following lemma.

\begin{lem}
Let $\pi=PG(2,q)$ with a distinguished line $l_{\infty}$.  Let $A$
be a subset of $l_{\infty}$. If the only functions for which $A$
contains the associated R\'{e}dei set are linear functions, then $A$
is primitive.
\end{lem}

\begin{rem}
We can also think of $\RR$ as the slope set (including infinity) or
the set of directions associated with $f$.
\end{rem}

\noindent Using Theorems \ref{T3} and \ref{T3b}  we get the
following.

\begin{cor}
Let $\pi=PG(2,q)$ with a distinguished line $l_{\infty}$. Suppose
$q=p^h$, $p$ prime, and let $t<h$ be maximal such that $t$ divides
$h$.  Let $A$ be a subset of $l_{\infty}$ contain the R\'{e}dei set
of the function $f$. If $|A|<\epsilon$ where
\begin{equation} \epsilon=\left\{\begin{array}{ll}
 \frac{q+3}{2} &  h=1\\
p^{h-t}+1 & \text{ otherwise.}
\end{array}\right.\end{equation}
then $f$ is a linear function.
\end{cor}

\noindent There is voluminous literature discussing R\'{e}dei sets
beginning with \cite{BrL}.  We mention also \cite{BBBSS},
\cite{Ball}, and \cite{szo2}.  The case where the set is a group is
discussed in \cite{BrL}.

\begin{thm}\label{L2} Let $\ell$ be a line in $\pi=PG(2,q)$ containing the
primitive set $A$.  Embed $\pi$ in $\Sigma=PG(3,q)$. Then $A$ is
primitive  in each plane of $\Sigma$ containing $\ell$.
\end{thm}
\begin{proof}
Let $\pi'$ be a plane of $\Sigma$ containing $\ell$.  Suppose by
way of contradiction that $S \subseteq \pi'\setminus \{\ell\}$ is
a set of $q$ points such that the line through any two points of
$S$ intersects $\ell$ in $A$ and that $S$ is not a subset of a
line. Select any point $P$ in $\Sigma\setminus\{\pi \cup \pi'\}$.
Let $\phi$ be the projection through $P$ mapping  $\pi'$ to $\pi$.
Then $\phi$ is a collineation fixing $\ell$.  Hence $\phi(S)$ is
not a subset of a line, yet any line on two points of $\phi(S)$
intersects $\ell$ within $A$. This contradicts the assumption that
$A$ is primitive in $\pi$.
\end{proof}

\noindent We generalize the property of being primitive to higher
dimensions as follows. Let $\Sigma=PG(k,q)$ let $\Pi$ be a
hyperplane of $\Sigma$ and let $E=\Sigma\setminus\Pi$ be the
associated affine space.  Let $S$ be a set of points of $E$. The
\textit{R\'{e}dei set} of $S$ with respect to $\Pi$, denoted by
$R_{\Pi}(S)$, is defined to be the set of all points of the form
$PQ\cap\Pi$, where $P,Q\in E$.

\begin{defn}\label{Dprim2} Let $\Sigma$, $\Pi$ and $S$ be as above with
$|S|=q^{k-1}$ and let $A$ be a subset of the points of $\Pi$. Then
$A$ is said to be \textit{primitive} if whenever
$R_{\Pi}(S)\subset A$, the set $S$ must be contained in a
hyperplane of $\Sigma$.
\end{defn}

\noindent The proof of the following is entirely similar to that
of Theorem \ref{L2}.

\begin{lem}
Let $\Pi$ be a hyperplane of $PG(k,q)$ containing the primitive
set $A$. Embed $PG(k,q)$ in $\Sigma=PG(k+1,q)$. Then $A$ is
primitive in every hyperplane of $\Sigma$ containing $\Pi$.
\end{lem}

\section{Linear Three-Dimensional MDS Codes}\label{l3d}

Let $C$ be a linear $(n,3,q)$-MDS code.  The associated generator
matrix $G$ is of rank $3$ over $\F=GF(q)$.  Again, let
$\Sigma=PG(3,q)$ be the underlying projective space of projective
dimension 3 as described in section \ref{LinBRS}. The $n$ columns
of $G$ give rise to a dual arc $\K$ in the plane at infinity,
denoted here by $\Pi$. Thus $\K=\{\ell_1,\ell_2,\ldots,\ell_n\}$
is a set of $n$ lines of $\Pi$ with no three collinear.  We want
to examine the extensions of $C$ to an MDS code.

A point of $\Pi$ lying on exactly two lines of $\K$ is called a
\textit{secant point}.  If the point lies on exactly one line of
$\K$ then it is called a \textit{tangent point}.  We note that
each line of $\K$ contains exactly $n-1$ secant points and
$q+1-(n-1)=q-n+2$ tangent points.

\begin{defn} Let $\K$ be a dual arc
in the projective plane $\pi$. Let $\ell$  be a line of $\K$ and
let $A$ denote the set of tangent points of $\K$ on $\ell$. Then
$\ell$ is said to be \textit{primitive with respect to $\K$} (or,
simply \textit{primitive}) if $A$ is a primitive set (in $\pi$).
\end{defn}

If $\K$ is a dual arc and $x$ is a line not in $\K$ such that
$\K\cup \{x\}$ is a dual arc then $x$ is said to \textit{extend}
$\K$, and $x$ is an \textit{extending line} of $\K$.

\begin{lem} \label{Lem1}
Let $\K$ be a dual $n$-arc in a projective plane of order $q$.  If
$n>q-\PP(q)+2$ then every line of $\K$ is primitive.
\end{lem}
\begin{proof}
Simply observe that each line of $\K$ is incident with precisely
$q-n+2<\PP(q)$ tangent points.
\end{proof}

Recall, (Definition \ref{dtrs}) that a set $S$ in $\Sigma-\Pi$ is
called a transversal set with respect to $\K$ if every line of
$\Sigma$ on a secant point of $\K$ intersects $S$ in at most one
point.

\begin{thm}\label{T5a}
Let $\K$ be a dual $n$-arc in $\pi=PG(2,q)$ with $\pi\subset
\Sigma=PG(3,q)$. Let $E=\Sigma-\pi$ be the associated affine
space. Suppose $S$ is  a transversal set of $\K$ with
$\abs{S}=q^2$. Assume there exist two primitive lines in $\K$.
Then $S$ is a hyperplane of $E$. Moreover, if $H$ is the
hyperplane of $\Sigma$ containing $S$, then $H \cap \pi$ is an
extending line of $\K$.
\end{thm}

\begin{proof}
Let $\ell$ and $\ell'$ be primitive with respect to $\K$. Let $\{\pi
_1,\pi_2,...\pi_q\}$ be the pencil of planes other than $\pi$
containing $\ell$.  Now $\abs{S}=q^2$.  Let $T= S \cap \pi_1$. We
claim $|T|=q$.  This follows since any of the $q$ lines of $\pi_1$
other than $\ell$ on a secant point of $\K$ can meet $T$ in at most
one point.  Thus, $\pi_1$ and indeed any of the planes $\pi_i$ meet
$S$ in at most, and hence in exactly, $q$ points.  By the
primitivity of $\ell$, $\pi_i\cap S$ must be an affine line $\ell_i$
in $\pi_i$, $1\leq i \leq q$.  Let $m_1,m_2,\ldots,m_q$ be the
corresponding projective lines, so that $m_i=\ell_i\cup P_i$ and
$m_i\cap\ell=P_i$, $1\leq i \leq q$, with $P_i$ on $\ell$.
\begin{enumerate}[(a)]
\item No two of the lines in
$\{m_1,m_2,\ldots,m_q\}$ are skew.  For suppose $m_1$ and $m_2$
are skew.  Then through every point $Q$ of $\pi$ off $\ell$, such
as a secant point, there is a line meeting $\ell_1$ and $\ell_2$.
This contradicts the fact that $S$ is a transversal set.  We
conclude that the lines $m_1,m_2,\ldots,m_q$ pass through a fixed
point $P$.
 \item  We claim that all of the lines $m_1,m_2,\ldots,m_q$ lie in
a plane.  To see this, fix a plane $\psi\neq \pi$ through $\ell'$.
As above $\psi$ meets $S$ in a line (through say $P'\in \ell'$).
Thus $\psi$ meets each of $\ell_1,\ell_2,\ldots,\ell_q$ in collinear
points $Q_1, Q_2,\ldots, Q_q$.  Therefore, the lines
$\ell_1,\ell_2,\ldots,\ell_q$ all lie in the plane $\gamma$
containing $P$ above and $Q_1, Q_2,\ldots, Q_q$. Thus, the points of
$S$ are the points of $\gamma$ in $E$. Moreover, since $S$ is a
transversal set, $\gamma$ meets $\pi$ in an extending line of $\K$.
\end{enumerate}
\end{proof}

\begin{rem}
Associated with $\ell$, $\ell'$ are points $P$, $P'$ as above.  We
note that $P\neq P'$ since $\ell\cap\ell'$ is a 2-fold point of $\K$
and $S$ is a transversal set of $\K$.  Also, $\ell'$ is not on $P$.
\end{rem}

\begin{cor}\label{cor1}
Let $\K$ be a dual $n$-arc in $\pi=PG(2,q)$, let $\Sigma=PG(3,q)$
and $E=\Sigma-\pi$ be the associated affine space. Let $S$ be  a
transversal set of $\K$ with $\abs{S}=q^2$. Assume that
$n>q+2-\PP(q)$. Then $S$ is a hyperplane of $E$.  Moreover, if $H$
is the hyperplane of $\Sigma$ containing $S$, then $H \cap \pi$ is
an extending line of $\K$.
\end{cor}

\begin{proof}
This follows immediately from  Lemma \ref{Lem1} and Theorem
\ref{T5a}.
\end{proof}

\begin{thm}\label{T7a}
Let $C$ be a linear  $(n,3,q)$-MDS code. Let $\K$ be a dual
$n$-arc giving rise to $C$ (as a BRS code).  If $\K$ contains two
primitive lines then any arbitrary $(n+1,3,q)$-MDS code extending
$C$ must be LE.
\end{thm}

\begin{proof}
As a BRS code, let $C$ be constructed within $\Sigma=PG(3,q)$.
Here $\,\pi$ is the hyperplane (= plane) at infinity containing
the dual $n$-arc $\K$, and $E=\Sigma-\pi$.

Let $C'$ be an $(n+1,3,q)$-MDS code extending $C$. Then $\,C'$
arises via a partition $\PP$=$C_1,C_2,\ldots C_q$ of $\,C$ where
each $C_i$ is an $(n,2,q)$-MDS code (Lemma \ref{Lptx}). The
partition $\PP$ corresponds to a partition
$\PP'=\{S_1,S_2,\ldots,S_q\}$ of $E$. Each $S_i$ is a set of
$q^{2}$ code words satisfying all conditions of $S$ in Theorem
\ref{T5a}. As such each $S_i$ is a set of $q^2$ points lying in a
plane $\pi_i$ of $\Sigma$, $1\leq i \leq q$.  Let $\pi_i\cap
\pi=\ell_i$.  As in the proof of Theorem \ref{T5a}, each $\ell_i$
is an extending line of $\K$.\\
\noindent Now any two planes $\pi_i$, $\pi_j$ of $\Sigma$ meet in a
line. Also, $\pi_i\cap E=S_i$ and $\pi_j\cap E=S_j$ are disjoint if
$i\neq j$, $1\leq i,j \leq q$.  It follows that
$\ell_1=\ell_2=\cdots =\ell_q=x$ and $x$ is an extending line of
$\K$.  Moreover, (by re-labelling the $(n+1)^{th}$ coordinate of
$C'$ if necessary) $C'$ is equivalent to the linear $(n+1,3,q)$-BRS
code associated with the dual arc $\K\cup\{x\}$.
\end{proof}

\noindent  With Lemma \ref{Lem1} and Theorem \ref{T7a} we get the
following.

\begin{thm}\label{xxa}
Let $C$ be a linear $(n,3,q)$-MDS code. If $n>q+2-\PP(q)$ then any
arbitrary  $(n+1,3,q)$-MDS code $C'$ extending $C$ must be LE.
\end{thm}

\begin{cor}\label{xx}
Let $C$ be a linear $(n,3,q)$-MDS code. Suppose $q=p^h$, $p$ prime,
and let $t<h$ be maximal such that $t$ divides $h$. If
$n>\beta$, where\\
 \[ \beta=\left\{
\begin{array}{ll}
\frac12(q+1) & \text{ if } q \text{ is prime,}\\
 q-p^{h-t}+1 & \text{ otherwise.}
\end{array}\right.\] then
any arbitrary  $(n+1,3,q)$-MDS code $C'$ extending $C$ must be LE.
\end{cor}

Let us discuss an easy application of this result.  The following
appears in \cite{Hi1} as Theorem 9.30.
\begin{thm}\label{Hir}
In $PG(2,q)$, $q\equiv 3$ (mod 4), there exist maximal $n$-arcs with
$n=\frac12(q+5)$. Equivalently, for such $q$ there exist linear
$(\frac12(q+5),3,q)$-MDS codes admitting no linear extensions.
\end{thm}

\noindent We can now state a stronger result for certain $q$ by
appealing to Corollary \ref{xx}.

\begin{lem}\label{non-extenders}
If $p\equiv 3$ (mod 4) is prime then  maximal
$(\frac{p+5}{2},3,p)$-MDS codes exist.
\end{lem}

\section{Higher Dimensions}

A dual $n$-arc $\K$ in $PG(k,q)$ is a collection
$\{\Pi=\Pi_1,\Pi_2,\ldots,\Pi_n\}$, $n\geq k+1$, of hyperplanes
such that no $k+1$ lie on a point, no $k$ lie on a line,..., no 3
lie on a $(k-2)$-flat, and no 2 lie on a $(k-1)$-flat.
Consequently, If we let $\Lambda_i=\Pi \cap \Pi_i$, $1<i\leq n$
then $\K'=\{\Lambda_1,\Lambda_2,\ldots,\Lambda_{n-1}\}$ is a dual
$(n-1)$-arc in $\Pi$. In this sense we say the remaining members
of $\K$ \textit{cut out} a dual $(n-1)$-arc in $\Pi$.

\begin{defn} Let $\K$ be a dual arc
in $\Pi=PG(k,q)$,$k\geq 2$. Let $\Lambda$  be a member of $\K$ and
let $B$ denote the set of $k$-fold points of $\K$. Then $\Lambda$
is said to be \textit{primitive} with respect to $\K$ (or, simply
primitive) if the point set $A=\Lambda\setminus B$ is a primitive
set.
\end{defn}

\begin{thm}\label{T5bb}
Let $\K$ be a dual $n$-arc in $\Pi=PG(k,q)$, $k\geq 2$.  Let
$\Sigma=PG(k+1,q)$ and let $S$ be a transversal set of $\K$ with
$\abs{S}=q^k$. If $\K$ contains two primitive members then $S$ is
a subset of a hyperplane $\HH$ of $\Sigma$. Moreover,  $\HH \cap
\Pi$ extends $\K$.
\end{thm}
\begin{proof}
Let $\Lambda$ be a primitive member of $\K$.  Let
$\{\Pi_1,\Pi_2,\ldots,\Pi_q\}$ be the pencil of hyperplanes of
$\Sigma$ other than $\Pi$ containing $\Lambda$.  As in the proof
of Theorem \ref{T5a} $\Pi_1$,  and indeed any of the hyperplanes
$\Pi_i$ meet $S$ in exactly $q^{k-1}$ points and (by the
primitivity of $\Lambda$) $\Pi_i\cap S$ is an affine $(k-1)$-flat
$\LL_i$ in $\Pi_i$, $1\leq i \leq q$.  Let $\M_1,\M_2,\ldots,\M_q$
be the corresponding projective $(k-1)$-spaces so that $\M_i\cap
\Lambda=\lambda_i$ and $\M_i=\LL_i\cup\lambda_i$, $1\leq i \leq
q$.

We claim the $\lambda_i$'s coincide.  Indeed, suppose the point
$P$ is in $\{\lambda_1\}\setminus\{\lambda_2\}$ and consider a
line $\ell$ in $\M_1$ where $\ell \cap \Pi=P$. In particular,
$\ell$ and $\M_2$ are disjoint.  Through each point $Q$ in $\Pi$
off $\Lambda$ there is a unique line meeting both $\ell$ and
$\M_2$ (necessarily in points of $S$). But then, since $S$ is a
transversal set, $Q$ can not be a $k$-fold point of $\K$.  Thus,
the $(k-1)$-spaces $\M_1,\M_2,\ldots,\M_q$ form a pencil on, say,
$\lambda$, where $\lambda=\lambda_1=\cdots=\lambda_q$.

We claim the $(k-1)$-spaces $\M_1,\M_2,\ldots,\M_q$ lie in a
hyperplane of $\Sigma$. Briefly, suppose $\M_1$, $\M_2$, and
$\M_3$ are not contained in a common hyperplane.  Let $\Pi_{12}$
be the unique hyperplane of $\Sigma$ containing $\M_1$ and $\M_2$.
Let $\Lambda'\neq \Lambda$ be a primitive member of $\K$.  Choose
a hyperplane $\Pi'$ other than $\Pi$ on $\Lambda'$. Denote by
$\tau_i$ the $(k-2)$-flat $\M_i\cap\Pi'$, $1\leq i \leq 3$. By
assumption $\tau_3$ is not contained in $\Pi_{12}$. By primitivity
there is a $(k-1)$-flat $H_1$ in $\Pi'$ containing $\tau_1,
\tau_2,$ and $\tau_3$. Let $H_2=\Pi'\cap\Pi_{12}$.  As $\tau_1$
and $\tau_2$ are (disjoint) $(k-2)$-flats, at most one
$(k-1)$-flat contains both. But then it follows that $H_1=H_2$ so
$\tau_3$ is contained in $\Pi_{12}$. The second conclusion of the
theorem is clear.
\end{proof}

\begin{thm}\label{T5bbb}
Let $\K$ be a dual $n$-arc in $\Pi=PG(k,q)$, $k\geq 2$.  If
$n>q-\PP(q)+k$ then every member of $\K$ is primitive.
\end{thm}
\begin{proof}
Our proof is inductive on $k$.  The case $k=2$ is given by Lemma
\ref{Lem1}.  Assume the result to hold in $PG(k-1,q)$. Let $\K$ be a
dual $n$-arc in $\Pi=PG(k,q), k>2$ with $n>q-\PP(q)+k$. Let $B$
denote the set of all $k$-fold points of $\K$. Choose $\lambda\in
\K$, let $A=\lambda\setminus B$ and suppose $S$ is a collection of
$q^{k-1}$ points of $\Pi-\lambda$ such that any line on a point of
$B$ intersects $S$ in at most one point (i.e. $R_{\lambda}(S)
\subset A$). The remaining members of $\K$ cut out a dual
$(n-1)$-arc $\K'$ in $\lambda=PG(k-1,q)$ so that in particular $S$
is a transversal set of $\K'$.  By the induction hypothesis, every
member of $\K'$ is primitive.  It follows (Theorem \ref{T5bb}) that
$S$ is contained in a hyperplane of $\Pi$ whence $\Lambda$ is
primitive.
\end{proof}

The proof of the following is entirely similar to that of Theorem
\ref{T7a}.

\begin{thm}\label{T8b}
Let $C$ be a linear $(n,k,q)$-MDS code, $k\geq 3$. Let $\K$ be a
dual $n$-arc giving rise to $C$ (as a BRS code).  If $\K$ contains
two primitive members then any arbitrary $(n+1,k,q)$-MDS code
extending $C$ must be LE.
\end{thm}

\noindent The last two theorems give the following.

\begin{thm}\label{C8b-a}
Let $C$ be a linear $(n,k,q)$-MDS code. If $n>q-\PP(q)+k$ then any
arbitrary  $(n+1,k,q)$-MDS code $C'$ extending $C$ must be LE.
\end{thm}

\begin{cor}\label{C8b}
Let $C$ be a linear $(n,k,q)$-MDS code. Suppose $q=p^h$, $p$ prime,
and let $t<h$ be maximal such that $t$ divides $h$. If $n>\beta$,
where \[ \beta=\left\{
\begin{array}{ll}
\frac12(q-3)+k & \text{ if $q$ is prime, and}\\
 q-p^{h-t}+k-1 & \text{ otherwise.}
\end{array}\right.\] then
any arbitrary  $(n+1,k,q)$-MDS code $C'$ extending $C$ must be LE.
\end{cor}

\section{Some Applications and Further Results} \label{lastsec}

We summarize some existing results and our corresponding
improvements.

\subsection{Applications for $q$ even}

\begin{thm}\label{T4}
Let $\K$ be an $n$-arc in $PG(2,q)$,  $n > \frac{q+2}{2}$ with $q$
even. Then $\K$ is contained in a unique maximal arc.
\end{thm}
\begin{proof}
See: \cite{szo}.
\end{proof}

The following Theorem is found in \cite{BBT}  with a proof using
results of algebraic geometry. For an inductive proof see
\cite{Wetl}.

\begin{thm}\label{TT5a}
Let $\K$ be an $n$-arc in $PG(k,q)$, $k\geq 2$ with $q$ even.  If
$n
>
 \frac{q}{2}+k-1 $. Then $\K$ is contained in a unique maximal
 arc.
\end{thm}

\noindent The above theorem together with Corollary \ref{C8b} gives
the following result.

\begin{thm}
Let $C$ be a linear $(n,k+1,q)$-MDS code over $GF(q)$, $q=2^h$. Let
$t<h$ be maximal such that $t$ divides $h$ and suppose $n>
q-2^{h-t}+k-1$.  Let $S$ be the collection of codes
consisting of $C$ and all extensions of $C$. Then \\
(i) all members of $S$ are LE; \\
(ii)  there is (up to equivalence) a
 unique  maximal code in $S$.\\
 \end{thm}

We now proceed to another application strengthening Theorem
\ref{hyper} below.

\begin{thm} \label{hyper}
Let $\K$ be an $n$-arc in $PG(2,q)$ with $n>q-\sqrt{q}+1$.  Then
$\K$ can be extended to a $(q+2)$-arc (a hyperoval) $\K'$ uniquely determined by $\K$.\\
Equivalently let $C$ be a linear $(n,3,q)$-MDS code over $GF(q)$
with $n>q-\sqrt{q}+1$. Let $S$ be the collection of codes consisting
of $C$ and all linear extensions of $C$. Then $S$ contains a
$(q+2,3,q)$-MDS code $C'$ and (up to equivalence) $C'$ is the only
maximal code in $S$.
\end{thm}
\begin{proof}
See \cite{Segr3}.
\end{proof}

\noindent Corollary \ref{C8b} can be used to strengthen Theorem
\ref{hyper} as follows.

\begin{thm}Let $C$ be a linear $(n,3,q)$-MDS code,  $n>q-\sqrt{q}+1$ with $q$ even. Then (up
to equivalence) there is exactly one maximal extension $C'$ of
$C$. Moreover $C'$ is an LE $(q+2,3,q)$-MDS code.
\end{thm}

\noindent Next we improve Theorem \ref{MC2} below.

\begin{thm} \label{MC2}
Let $\K$ be a $(q+1)$-arc in $PG(k,q)$, $q$ even, with either (a)
$k=3,4$ ; or
(b) $k\geq 5$ and $q \geq (k-2)^3$.  Then $\K$ is complete.\\
Equivalently let $C$ be a linear $(n,k+1,q)$-MDS code satisfying (a)
or (b). Then $C$ can not be extended to a linear $(q+2,k+1,q)$-MDS
code.
\end{thm}
\begin{proof}
For (a) see \cite{Cas} and for (b) see \cite{BTB}
\end{proof}

\noindent Corollary \ref{C8b} strengthens the above to the
following:

\begin{thm} Let
$C$ be a linear $(q+1,k+1,q)$-MDS code with either (a) $k=3,4$ ;
or (b) $k\geq 5$ and $q \geq (k-2)^3$.  Then $C$ is maximal.
\end{thm}

\begin{thm}\label{NRC}
Let $\K$ be an $n$-arc in $PG(k,q)$ $k\geq 4$ and $n\geq
q-\sqrt[3]{q}+k$.  Then $\K$ lies in a normal rational curve
uniquely determined by $\K$.\\
Equivalently:\\
Let $C$ be a linear $(n,k+1,q)$-MDS code over $GF(q)$ with $k\geq
4$ and $n\geq q-\sqrt[3]{q}+k$. Let $S$ be the collection of codes
consisting of $C$ and all linear extensions of $C$. Then up to
equivalence there is an unique maximal code $C'$ in $S$.  Moreover
$C'$ is (equivalent to) a GRS-code.
\end{thm}
\begin{proof}
See \cite{BTB}.
\end{proof}

\noindent Using Corollary \ref{C8b} we have the following
improvement to Theorem \ref{NRC}.

\begin{thm}
Let $C$ be a linear $(n,k+1,q)$-MDS code over $GF(q)$ with $k\geq
4$ and $n\geq q-\sqrt[3]{q}+k$.  Then up to equivalence there is
an unique maximal extension $C'$ of $C$.  Moreover $C'$ is
(equivalent to) a GRS-code.
\end{thm}

\subsection{Applications for $q$ odd}

\begin{thm}\label{TT4}
Let $\K$ be an $n$-arc in $PG(2,q)$ with $n > \frac23(q+2)$ and
$q$ odd. Then $\K$ is contained in a unique maximal arc.
\end{thm}
\begin{proof}
See: \cite{szo}.
\end{proof}

The following can be found in \cite{BBT}.

\begin{thm}\label{T5b}
Let $\K$ be an $n$-arc in $PG(k,q)$, $k \geq 2$ with $q$ odd. If
$n>\frac{2}{3}(q-1)+k$. Then $\K$ is contained in a unique
maximal arc.\\
 Equivalently:\\
Let $C$ be a linear $(n,k+1,q)$-MDS code over $GF(q)$ with $n >
\frac{2}{3}(q-1)+k$ and let $S$ be the collection of codes
consisting of $C$ and all linear extensions of $C$.  Then (up to
equivalence) there is a unique maximal code in $S$.
\end{thm}

\noindent Corollary \ref{C8b} strengthens the above to the following
in the case that $q$ is prime.

\begin{thm}
For $q$ an odd prime, let $C$ be a linear $(n,k+1,q)$-MDS code over
$GF(q)$ with $n> \frac{2}{3}(q-1)+k$.  Let $S$ be the collection
of codes consisting of $C$ and all extensions of $C$. Then:\\
(i) All members of $S$ are LE; and\\
(ii)  There is (up to equivalence) a
 unique  maximal code in $S$.\\
 \end{thm}

\noindent If $q$ odd is not prime then $q \geq 9$, so $n$ integral
gives
\[ n>q-\sqrt{q}+k-1 \Rightarrow n>\frac23(q-1)+k.\]
Hence, Theorem \ref{T5b} and Corollary \ref{C8b} give the following.

\begin{thm}
For $q$ odd not a prime, let $C$ be a linear $(n,k+1,q)$-MDS code
over $GF(q)$.  Let $S$ be the set of codes consisting of $C$ and
all extensions of $C$.  If $n> q-\sqrt{q}+k-1$ then:\\
(i) All members of $S$ are LE; and\\
(ii)  There is (up to equivalence) a
 unique  maximal code in $S$.
 \end{thm}

\noindent The previous result can be considerably improved in most
cases by observing that if $q=p^h$ and $t<h$ with $p^t>3$ then
\[ n>q-p^{h-t}+k-1 \Rightarrow n>\frac23(q-1)+k.\]

\noindent Hence, Theorem \ref{T5b} and Corollary \ref{C8b} give the
following.

\begin{thm}
Let $C$ be a linear $(n,k+1,q)$-MDS code over $GF(q)$, where $q=p^h$
is odd. Let $t<h$ be maximal such that $t$ divides $h$. Let $S$ be
the set of codes consisting of $C$ and all extensions of $C$. If
$p^t>3$ and $n> q-p^{h-t}+k-1$ then:\\
(i) All members of $S$ are LE; and\\
(ii)  There is (up to equivalence) a
 unique  maximal code in $S$.
 \end{thm}

\noindent We now want to improve Theorem \ref{abc} below.

\begin{thm} \label{abc}
Let $\K$ be a dual $(q+1)$-arc in $PG(k,q)$, $q$ odd,  with either\\
(a) $k=2,3,$ or 4; or (b) $q>(4k-23/4)^2$. Then $\K$ is complete.\\
Equivalently:\\
Let $C$ be a linear $(q+1,k+1,q)$-MDS code satisfying (a) or (b)
above, then $C$ can not be extended to a linear $(q+2,k,q)$-MDS
code.
\end{thm}
\begin{proof}
For (a) see \cite{Segr} and for (b) see \cite{JAT1,JAT2}
\end{proof}

\noindent Corollary \ref{C8b}, combined with Theorem \ref{abc} now
gives the following improvement of Theorem \ref{abc}.

\begin{thm} Let $C$ be a linear $(q+1,k+1,q)$-MDS code, $q$ odd, with
(a) $k=2,3,$ or 4; or (b) $q>(4k-23/4)^2$.
 Then $C$ is maximal.
\end{thm}

\begin{thm} \label{GRS2}
Let $\K$ be a $n$-arc in $PG(k,q)$, $q$ odd,  with
$n>q-\frac14\sqrt{q}+k-\frac74$ and either (a) $ k=2,3$; or (b)
$q>(4k-\frac{23}{4})^2$.  Then $\K$ is contained in an unique
(necessarily maximal) $(q+1)$-arc $\K'$.  Moreover $\K'$ is the
point set of a normal rational curve.\\
Equivalently:\\
Let $C$ be a linear $(n,k+1,q)$-MDS code over $GF(q)$  with
$n>q-\frac14\sqrt{q}+k-\frac74$ and either (a) k=2,3; or (b)
$q>(4k-\frac{23}{4})^2$. Let $S$ be the collection of codes
consisting of $C$ and all linear extensions of $C$. Then $S$
contains an unique maximal code $C'$, moreover $C'$ is (equivalent
to) a GRS-code.
\end{thm}
\begin{proof}
If $k=2$ or 3 then every $q+1$ arc is the point set of a normal
rational curve (\cite{Segr,Segr3}). For any $k$, if
$q>(4k-\frac{23}{4})^2$ then every $q+1$ arc is the point set of a
normal rational curve (\cite{JAT1,JAT2}).  For any $k$, if
$n>q-\frac14\sqrt{q}+k-\frac74$ than $\K$ is contained in an
unique normal rational curve (\cite{JAT1,JAT2}).
\end{proof}

\noindent Corollary \ref{C8b} strengthens Theorem \ref{GRS2} as
follows.

\begin{thm}
Let $C$ be a linear $(n,k+1,q)$-MDS code over $GF(q)$  with $q$
odd, $n>q-\frac14\sqrt{q}+k-\frac74$, and either (a) k=2,3; or (b)
$q>(4k-\frac{23}{4})^2$.  Then (up to equivalence) there is an
unique maximal code $C'$ extending $C$.  Moreover, $C'$ is
(equivalent to) a GRS-code.
\end{thm}


\end{document}